\documentclass[11pt]{article}

\usepackage{amsfonts}
\usepackage{amsmath,amssymb,latexsym,a4}

\RequirePackage{hyperref}
\RequirePackage{hyperref}
		\usepackage{pgfplots}
		\usetikzlibrary{plotmarks}
		\usetikzlibrary{external}
		\pgfplotsset{compat=1.3}
		\newlength\figurewidth 
		\newlength\figureheight
\usepackage{epsf}
\usepackage{color}
\usepackage{epsfig}

\usepackage{graphicx}
\usepackage{amsmath}

\newcommand{\e}{{\rm e}}

\usepackage{graphicx}
\usepackage{mathdots}
\usepackage{array}

		\usepackage{pgfplots}
		\usetikzlibrary{plotmarks}
		\usetikzlibrary{external}
		\pgfplotsset{compat=1.3}

\makeatletter
\def\Ddots{\mathinner{\mkern1mu\raise\p@
\vbox{\kern7\p@\hbox{.}}\mkern2mu
\raise3\p@\hbox{.}\mkern2mu\raise5\p@\hbox{.}\mkern1mu}}
\makeatother
\usepackage{booktabs}

\usepackage{amsmath,amssymb}

\def\en{\text{\sc{e}-}}

\def\mdim{N} 

\def\e{e}
\def\id{I}

\def\C{\mathbb{C}}

\def\beq{\begin{equation}}
\def\eeq{\end{equation}}

\begin{document}

\title{An improved algorithm to compute the exponential of a matrix}

	\author{Philipp Bader\footnote{Departament de Matemàtiques, Universitat Jaume I, 12071 Castellón, Spain. {\tt  bader@uji.es}},
		\ \ Sergio Blanes\footnote{Instituto de Matem\'atica Multidisciplinar,
  Universitat Polit\`ecnica de Val\`{e}ncia, E-46022  Valencia, Spain.
 {\tt serblaza@imm.upv.es}}, 
		\ \ Fernando Casas\footnote{IMAC and Departament de Matemàtiques, Universitat Jaume I, 12071 Castellón, Spain. {\tt  fernando.casas@uji.es}}}

\maketitle

\begin{abstract}
In this work, we present a new way to compute the Taylor polynomial of the matrix exponential which reduces the number of matrix multiplications in comparison with the de-facto standard Patterson--Stockmeyer method.
This reduction is sufficient to make the method superior in performance to Pad\'e approximants by 10-30\% over a range of values for the matrix norms and thus we propose its replacement in standard software kits.
Numerical experiments show the performance of the method and illustrate its stability.
\end{abstract}



\vspace*{.5cm}
\section{Introduction}

The search for efficient algorithms to compute the exponential of an square matrix has a tremendous interest, given the wide range of its applications in many branches of science. Its
importance is clearly revealed by the great impact achieved by various reference reviews devoted to the subject, e.g.  \cite{higham05tsa,moler03ndw,najfeld95dot,sidje98eas}. 
The problem has been tackled of course by many authors and a wide variety of techniques have been proposed. Among them, the scaling and squaring procedure is perhaps the most
popular when the dimension of the corresponding matrix runs well into the hundreds. As a matter of fact, this technique is incorporated in popular computing packages such as
{\sc Matlab} (\texttt{expm}) and {\sc Mathematica} (\texttt{MatrixExp}) in combination with Pad\'e approximants \cite{almohy09ans,higham05tsa,higham10cma}. 

Specifically, for a given matrix $A \in\C^{\mdim\times \mdim}$, the scaling and squaring technique is based on the key property
\begin{equation}  \label{eses}
	\e^A = \left( \e^{A/2^s}  \right)^{2^s} 
	, \quad s\in\mathbb{N}.
\end{equation}
The exponential $\e^{A/2^s}$ is then replaced by a rational approximation, namely a $2m$th-order diagonal Pad\'e approximant, $r_{2m}(A/2^s)$.
The optimal choice of both parameters, $s$ and $2m$, depends on the estimation of a bound of $\|A\|$ and the desired accuracy. 
%
Diagonal Pad\'e approximants have the form 
\begin{equation}\label{eq.24}
		r_{2m}(A )=\frac{p_{m}(A )}{p_{m}(-A )},
\end{equation}
where the polynomial $p_m(x)$ is given by
\[
   p_{m}(x)= \sum_{j=0}^m \frac{(2m - j)! m!}{(2m)! (m-j)!} \frac{x^j}{j!}
\]
so that  $r_{2m}(A) = \e^{A }+\mathcal{O}(A^{2m+1})$ and the evaluation of $p_{m}(A ), p_{m}(-A )$ is carried out with a reduced number of matrix products. Finally,  
the cost of the inverse is taken as $4/3$ the cost of one product: it requires one $LU$ factorization at the cost of $1/3$ products plus $N$ solutions of upper and lower triangular systems by forward and backward substitution at the cost of one product. In particular, $r_{10}(A)$ and $r_{26}(A)$ are considered among the optimal choices (with respect to accuracy and computational cost) in this approach when high accuracy is desired and the norm $\|A\|$ is relatively large. The corresponding algorithm for $r_{10}(A)$ is given by  
 \cite{higham05tsa}
\begin{equation}\label{eq.pade10}
\begin{aligned}
   & 	u_{5}  =  A \big( b_5 A_4 +  b_3 A_2 + b_1 \id \big), \\
  &  v_{5} =  b_4 A_4 + b_2 A_2 + b_0 \id,  \\
   &  (-u_{5}+v_{5}) \, r_{10}(A) = u_{5}+v_{5},
\end{aligned}   
\end{equation}
whereas for   $r_{26}(A)$ one has
\begin{equation}\label{eq.pade26}
\begin{aligned}
   &  	u_{13}  =  A \Big( A_6 \big( b_{13}A_6 + b_{11}A_4 + b_9A_2 \big) + b_7A_6 + b_5A_4 + b_3A_2 + b_1 \id \Big), \\
   &	v_{13} =  A_6 \big( b_{12}A_6 + b_{10}A_4 + b_8A_2 \big)  + b_6A_6 + b_4A_4 + b_2A_2 + b_0\id, \\
   &   	(-u_{13}+v_{13}) r_{26}(A) = u_{13}+v_{13}.
 \end{aligned}
 \end{equation}  
Here $A_2=A^2$, $A_4=A_2^2$ and $A_6=A_2A_4$.
Written in this form, it is clear that only three and six matrix multiplications and one inversion are required to obtain approximations of order 10 and 26, respectively. 
These two methods are used in fact in the current implementation of the function  \verb"expm"  in {\sc Matlab}.


The choice of the optimal order of the approximation for a given value of $\|A\|$ is based on the control of the backward error. To be more specific, given a prescribed accuracy
$u$, the parameter $m$ is selected according to $\theta_{2m}(u)$, the largest value of $\|A\|$ such that the Pad\'e scheme $r_{2m}$ has relative backward-error smaller than $u$, i.e.,
\begin{equation}\label{eq.theta}
	\forall A, \ \  \|A\|\leq \theta_{2m}(u) \ :\  r_{2m}(A)=e^{A+E}, \quad \text{s.t.}\  \frac{\|E\|}{\|A\|}\leq u.
\end{equation}
The values of $\theta_{2m}$ for $u=2^{-24}$ and $u=2^{-53}$, corresponding to single and double precision, respectively, are collected in 
Table~\ref{tab.theta}. They have been obtained according to 
 \cite[Theorem 2.1]{higham05tsa} by defining the majorant of 
 \[
    \rho(\theta)=\log(e^{-\theta}r_{2m})=\sum_{i=2m}^\infty c_i \theta^i
\]
to $f(\theta)=\sum_{i=2m}^\infty |c_i|\theta^i$ and solving $f(\theta))/\theta=u$, where the series is truncated after 150 terms.

\begin{table}\centering\footnotesize
\caption{\label{tab.theta} Values of $\theta_{2m}$ for diagonal Pad\'e of order $2m$ with the minimum number of products $\pi$ for single and double precision.
}
\newcolumntype{H}{@{}>{\lrbox0}l<{\endlrbox}}
\newcolumntype{D}{>{$}r<{$}}
\begin{tabular}{DDDDHDHDHDHHHD}
\toprule
\pi: & 0 & 1 & 2 & 3&3&4&4&5&5&5&6&6&6\\
	 2m:  &2					&  			4
	&		  6 &  		  8&		   10&  		  12&       14&   16&   18&  20&  22&  24&  26 
	\\ \midrule	 
	u \leq2^{-24}&
	8.46\en4&
	8.09\en2&
	4.26\en1 &
	1.05&
	1.88&
	2.85&
	3.93&
	5.06&
	6.25&
	7.47&
	8.71&
	9.97&
	11.2\\
u\leq2^{-53}&
				3.65\en8&				5.32\en4&	
				1.50\en2&8.54\en2&2.54\en1&5.41\en1&9.50\en1&1.47&2.10&2.81&3.60&4.46& 5.37 
				\\
\bottomrule
\end{tabular}
\end{table}


Another technique that has recently been proposed consists of using, instead of a diagonal Pad\'e approximant, the Taylor polynomial of degree $m$, 
$T_m(A) = \sum_{k=0}^mA^k/k!$, evaluated according
with the Paterson--Stockmeyer procedure (see \cite{higham08fom,higham10cma,paterson73otn}). In that way the number of matrix products is reduced and the overall
performance is improved for matrices of small norm, although it is less efficient for matrices with large norm \cite{higham10cma,ruiz16hpc,sastre11ame,sastre14aae,sastre15nsq}.

If the Paterson--Stockmeyer (PS) technique is carried out to compute $T_m(A)$ in a Horner-like fashion, the maximal attainable degree is $m=(k+1)^2$ by using $2k$ matrix
products. The optimal choices for most cases then correspond to $k=2$ (four products) and $k=3$ (six products), i.e, to degree $m=9$ and $m=16$, respectively. The
corresponding polynomials are then computed as follows:
\begin{equation}   \label{poli9}
\begin{aligned}
   & T_{9}(A)=\sum_{i=0}^9 c_i A^i = f_0+ \big(f_1+(f_2+c_{9} A_3) A_3 \big) A_3,  \\
   & T_{16}(A)=\sum_{i=0}^{16} c_i A^i = g_0+ \big(g_1+(g_2+(g_3+c_{16} A_4) A_4) A_4 \big) A_4,
\end{aligned}
\end{equation}
where $c_i=1/i!$ and
\[
\begin{aligned}
  & f_i=\sum_{k=0}^2 c_{3i+k} A_k,  \qquad i=0,1,2   \\
 &  g_i=\sum_{k=0}^3 c_{4i+k} A_k,  \qquad i=0,1,2,3.
\end{aligned}
\]
Here, as before,  $A_2=A^2$, $A_3=A^2 A$, $A_4=A_2 A_2$. In Table~\ref{tab.theta.new} we collect the values for the corresponding thresholds $\theta_m(u)$ 
needed to select the best scheme for a given accuracy.

The main goal of this work is to show that it is indeed possible to organize the computation of the Taylor polynomial of the matrix exponential function  in a more efficient way than
the Paterson--Stockmeyer technique, so that  with the same
number of matrix products one can construct a polynomial of higher degree. In this way, a method for computing $\e^A$ based on scaling and squaring together with the new
implementation of the Taylor polynomial is proposed which is more efficient than the same procedure based on Pad\'e approximants for a wide range of values of $\|A\|$.


\begin{table}\centering\footnotesize
\caption{\label{tab.theta.new} Values of $\theta_{m}$ for Taylor of degree $m$ the with minimum number of products $\pi$ for single and double precision.
{We have included degree 24 to illustrate that the gain is marginal over scaling and squaring for double precision and negative for single precision. The $\theta_m$ values should at least double per extra product to be favorable over scaling and squaring.}
}
\newcolumntype{H}{@{}>{\lrbox0}l<{\endlrbox}}
\newcolumntype{D}{>{$}r<{$}}
\begin{tabular}{DDDDHDHDHDHD}
\toprule
	\pi & 0 & 1 & 2 & 3 & 3 &4& 4 & 5&5 &6& {6} \\
	m	&1					&  			2
	&		  4 &  		  6&	 8& 	   9&    12&  		  16&      18&    20& 24
	\\ \midrule
u \leq2^{-24}&1.19\en7 & 5.98\en4&
				5.12\en2 & 2.50\en-1& 5.80\en1&7.80\en1&1.46&2.48&3.01&3.55&4.65\\
 u \leq2^{-53}& 2.22\en16 & 2.58\en8&
				3.40\en4&9.07\en3&4.99\en2&8.96\en2&2.99\en1&7.80\en1&1.09&1.44&2.22 
				\\
\bottomrule
\end{tabular}
\end{table}

\section{A generalized recursive algorithm}

Clearly, the most economic way to construct polynomials of degree $2^k$ is by applying the following sequence, which involves only $k$ products:
\begin{align}
A_1 &= A \nonumber  \\
A_2 &= (\delta_1I+x_1A_1)(\delta_2 I+x_2 A_1)  \nonumber \\
A_4 &= (\delta_3I+x_3 A_1+x_4A_2)(\delta_4I+x_5A_1+x_6A_2)  \label{New_algorithm}
\\
A_8 &= (\delta_5I+x_7A_1+x_8A_2+x_9A_4)(\delta_6I+x_{10}A_1+x_{11}A_2+x_{12}A_4),\nonumber \\
 & \vdots  \nonumber 
\end{align}
with $\delta_i=0,1$. These polynomials are then linearly combined to form
\[
T_{2^k} = y_0I+y_1A_1+y_2A_2+y_3A_4+y_4A_8+\cdots +y_{k+1}A_{2^k}.
\]
Notice that the indices in $A$, $A_{2^k}$, are chosen to indicate the highest attainable power, i.e., $A_{2^k}=\mathcal{O}(A^{2^k})$.
A simple counting tells us that with $k$ products one has $(k+1)^2+1$ {parameters} to construct a polynomial of degree $2^k$ that contains $2^k+1$ {coefficients}. It is then clear 
that the number of coefficients grows faster than the number of parameters, so that this procedure cannot be used to obtain high degree polynomials, as already noticed in \cite{{paterson73otn}}. Even worse, in general, not all parameters are independent and this simple estimate does not suffice to guarantee the existence of solutions with real coefficients. 

Nevertheless, this procedure can be modified in such a way that additional parameters are introduced, at the price, of course, of including some extra products. In particular, we could
include new terms of the form 
\[
  (\gamma_1I+z_1A_1)(\gamma_2I+z_2A_1+z_3A_2),
\]
not only in the previous $A_k$, $k > 1$, but also in $T_{2^k}$, which would allow us to introduce a cubic term and an additional parameter.

Although the Paterson-Stockmeyer (PS) technique is arguably the most efficient procedure to evaluate a \emph{general} polynomial, there are relevant classes of
 polynomials for which the PS rule involves more products than strictly necessary. To illustrate this feature, let us consider the evaluation of 
\begin{equation}   \label{pol.w}
 \Psi(k,A)=I+A+\cdots+A^{k-1},
\end{equation}
a problem addressed
 in \cite{westreich89etm}. Polynomial (\ref{pol.w}) appears in connection with the integral of the state transition matrix and the analysis of multirate sampled data systems. In
\cite{westreich89etm} it is shown that
 with three matrix products one can evaluate $\Psi(7,A)$ (as with the PS rule), whereas with four products it is possible to compute $\Psi(11,A)$
(one degree higher than using the PS rule). In general, the savings with respect to the PS technique grow with the degree $k$. The procedure was further improved and analysed in \cite{li92afa}, where the following conjecture was posed: the minimum number of products to evaluate  $\Psi(k,A)$ is $2\left\lfloor \log_2 k \right\rfloor -2 + i_{j-1}$ where $N=(i_j,i_{j-1},\ldots,i_1,i_0)_2$ (written in binary), i.e., $i_{j-1}$ is the second most significant bit.

This conjecture is not true in general, as is illustrated by the following algorithm of type \eqref{New_algorithm}, 
that allows one to compute  $\Psi(9,A)$ by using only three matrix products: 
\begin{align}
A_2 &= A^2, \qquad B=x_1I+x_2A+x_3A_2  \nonumber \\
A_4 &= x_4I+x_5A+B^2  \label{Algorithm83}
\\
A_8 &= x_6A_2+A_4^2  \nonumber  \\
\Psi(9,A) & = x_7I+x_8A+x_9A_2+A_8,  \nonumber
\end{align}
with
\[
\begin{array}{ccccc}
 x_1=  \displaystyle \frac{-5 + 6 \sqrt{7}}{32},  & x_2  \displaystyle =-\frac14,  & x_3  =-1, & 
x_4=   \displaystyle  \frac{3 (169 + 20 \sqrt{7})}{1024}, & \\
 x_5= \displaystyle  \frac{3(5 + 2\sqrt{7})}{64}, & 
x_6= \displaystyle \frac{1695}{4096},  & x_7= \displaystyle  \frac{267}{512},  & x_8= \displaystyle \frac{21}{64},  & x_9=  \displaystyle \frac{3 \sqrt{7}}4.
\end{array}
\]
Notice that, since the coefficients of the algorithm must satisfy a system of nonlinear equations, they are irrational numbers.

Although by following this approach it is not possible to achieve degree 16 with four products, there are other polynomials of degree 16 that can indeed be computed with 
only four products. This is the case, in particular, of the 
 truncated Taylor expansion of the function $\cos(A)$:
\[
\displaystyle T_{16}=\sum_{i=0}^8 \frac{(-1)^iA^{2i}}{(2i)!}=\cos(A)+{\cal O}(A^{17}).
\]
Taking $B=A^2$ we obtain a polynomial of degree 8 in $B$ that can be evaluated with three additional products in a similar way as in the computation of $\Psi(9,A)$, 
but with different coefficients.

\section{A new algorithm to compute the matrix exponential}

Algorithm  (\ref{New_algorithm}) can be conveniently modified along the lines exposed in the previous section to compute the 
truncated Taylor expansion of the matrix exponential function
\begin{equation}   \label{tay.exp}
\displaystyle T_{n}(A)=\sum_{i=0}^n \frac{A^{i}}{i!}=\e^A+{\cal O}(A^{n+1})
\end{equation}
for different values of $n$ using the minimum number of products. In practice,
 we proceed in the reverse order: given a number $k$, we find a convenient (modified) sequence of type  (\ref{New_algorithm}) that allows one to construct the highest degree polynomial $T_n(A)$ using only $k$ matrix products.
 The coefficients in the sequence satisfy a relatively large system of algebraic nonlinear equations. Here several possibilities may occur:  (i) the system has no solution; (ii) there is
 a finite number of real and/or complex solutions, or (iii) there are families of solutions depending on parameters. For our purposes, we only need one solution with real coefficients. 
 In addition, if there are several solutions we take a solution with small coefficients to avoid large round off errors due to products of large and small numbers.

With $k=0,1,2$ products we can evaluate  $T_n$ for $n=1,2,4$, in a similar way as the PS rule, whereas for $k=3,4,5$ and 6 products the situation is detailed next.

\

\paragraph{$k=3$ products}
In this case, with the PS rule only $T_6$ can be determined, whereas the following algorithm allows one to evaluate  $T_8$:
\begin{equation}  \label{Algorithm83eA}
\begin{aligned}
   A_2 & = A^2\\
   A_4 & = A_2(x_1 A+x_2 A_2)\\
   A_8 & = (x_3 A_2+A_4)(x_4I+x_5A+x_6 A_2+x_7A_4)  \\
   T_8(A)  & = y_0 I + y_1 A +y_{2} A_2+ A_8. 
\end{aligned}  
\end{equation}
With this sequence we get two families of solutions depending on a free parameter, $x_3$, which is chosen to (approximately) minimize the 1-norm of the vector of parameters. The reasoning behind this approach is to achieve numerical stability by avoiding the multiplication of high powers of $A$ with large coefficients, in a similar vein as the Horner procedure. 
The coefficients in (\ref{Algorithm83eA}) are given by
\[
\begin{array}{llll}
x_1=\displaystyle x_3\frac{ 1 + \sqrt{177}}{88},& 
x_2= \displaystyle \frac{1 + \sqrt{177}}{352}{x_3},&
x_4= \displaystyle \frac{-271 + 29\sqrt{177}}{315 x_3}, \\
x_5= \displaystyle \frac{11 (-1 + \sqrt{177})}{1260 x_3},&
x_6=  \displaystyle \frac{11 (-9 + \sqrt{177})}{5040 x_3},&
x_7=  \displaystyle \frac{89 - \sqrt{177}}{5040 x_3^2},&\\
y_0=1,& 
y_1= 1,& 
y_2 = \displaystyle \frac{857 - 58\sqrt{177}}{630}, \\
x_3 = 2/3.
\end{array}
\]
Perhaps surprisingly, $T_7(A)$ requires at least 4 products, so $T_8$ may be considered as a singular polynomial.


\paragraph{$k=4$ products}
Although polynomials up to degree 16 can in principle be constructed, our analysis suggests that the Taylor polynomial (\ref{tay.exp}) corresponding to $\exp(A)$ does not
belong to that family. To determine the highest degree one can achieve with $k=4$, we proceed as follows.
 We take  $T_n(A)$ for a given value of $n$ and decompose it as a product of two polynomials of lower degree plus a lower degree polynomial (that will be used to evaluate the higher degree polynomials). 
The highest value we have managed to reach is $n=12$. It is important to note that our ansatz gives many different ways to write the sought polynomial, one of which is given by
			\begin{align*}
A_2 &= A^2  \\
			A_3 &= A_2A, \qquad\qquad B=x_1I+x_2A_1+x_3A_2+x_4A_3\\
			A_6 &= x_5I+x_6A_1+x_7A_2 - B^2\\
			A_{12} &= (x_8A_2 + x_9A_3 + A_6)A_6\\
		T_{12} &{= y_0 I + y_1 A_1 + y_2 A_2 + y_3 A_3 + A_{12}} .
			\end{align*}
	This approach, however, leads to a much less stable numerical behaviour.
	We propose instead to use the following sequence, which has been empirically shown to be comparable in stability to Pad\'e and Horner methods.
\begin{equation}  \label{poly12}
	\begin{aligned}
		A_2&=A^2,\\
		A_3&=A_2A,\\
		B_1 &= a_{0,1}I+a_{1,1}A+a_{2,1}A_2+a_{3,1}A_3,\\
		B_2 &= a_{0,2}I+a_{1,2}A+a_{2,2}A_2+a_{3,2}A_3,\\
		B_3 &= a_{0,3}I+a_{1,3}A+a_{2,3}A_2+a_{3,3}A_3,\\
		B_4 &= a_{0,4}I+a_{1,4}A+a_{2,4}A_2+a_{3,4}A_3,\\
		A_6 &= B_3 + B_4^2 \\
		T_{12}(A) & = B_1 + (B_2 + A_6)A_6.
	\end{aligned}
	\end{equation}
	The algebraic expressions which minimize the 1-norm of the vector formed by the parameters do not give additional insight and hence we show their numeric values only:
	\[\begin{array}{lrlrlrlrlrl}
	a_{0,1} =& 9.0198\cdot 10^{-16},&
	a_{0,2} =& 5.31597895759871264183,\\
	a_{0,3} =& 0.18188869982170434744,&
	a_{0,4} =&-2.0861320\cdot 10^{-13},\\
	a_{1,1} =& 0.46932117595418237389,&
	a_{1,2} =& 1.19926790417132231573,\\
	a_{1,3} =& 0.05502798439925399070,&
	a_{1,4} =&-0.13181061013830184015,\\
	a_{2,1} =&-0.20099424927047284052,&
	a_{2,2} =& 0.01179296240992997031,\\
	a_{2,3} =& 0.09351590770535414968,&
	a_{2,4} =&-0.02027855540589259079,\\
	a_{3,1} =&-0.04623946134063071740,&
	a_{3,2} =& 0.01108844528519167989,\\
	a_{3,3} =& 0.00610700528898058230,&
	a_{3,4} =&-0.00675951846863086359.
	\end{array}
	\]

\paragraph{$k=5$ products}
With 5 products, $n=18$ is the highest value we have been able to achieve. We write $T_{18}$ as the product of two polynomials of degree 9, 
that are further decomposed into polynomials of lower degree. The polynomial is evaluated through the following sequence:
\begin{equation}   \label{poly.18}
\begin{aligned}
A_2 &= A^2, \quad A_3 = A_2A, \quad A_6 = A_3^2,\\
B_{1} &= a_{0,1}I + a_{1,1}A + a_{2,1}A_2 + a_{3,1}A_3,\\
B_{2} &= b_{0,1} I + b_{1,1}A + b_{2,1}A_2 + b_{3,1}A_3 + b_{6,1}A_6, \\
B_{3} &= b_{0,2} I + b_{1,2}A + b_{2,2}A_2 + b_{3,2}A_3 + b_{6,2}A_6,\\
B_{4} &= b_{0,3} I + b_{1,3}A + b_{2,3}A_2 + b_{3,3}A_3 + b_{6,3}A_6,\\
B_{5} &= b_{0,4} I+  b_{1,4}A + b_{2,4}A_2 + b_{3,4}A_3 + b_{6,4}A_6,\\
A_{9}  &= B_{1}B_{5} + B_{4},\\
T_{18}(A) &= B_{2} + (B_{3} + A_9)A_9,
\end{aligned}								
\end{equation}
with coefficients 
\[
\begin{array}{>{}lrlr}
a_{0,1} =& 0, &
a_{1,1} =& -0.100365581030144620,\\
a_{2,1} =&-0.0080292464824115696,&
a_{3,1} =&-0.0008921384980457299,\\
b_{0,1} =& 0,&
b_{1,1} =&0.39784974949964507614,\\
b_{2,1} =&1.36783778460411719922,&
b_{3,1} =&0.49828962252538267755,\\
b_{6,1} =&-0.0006378981945947233,&
b_{0,2} =&-10.967639605296206259,\\
b_{1,2} =&1.68015813878906197182,&
b_{2,2} =&0.05717798464788655127,\\
b_{3,2} =&-0.0069821012248805208,&
b_{6,2} =&0.00003349750170860705,\\
b_{0,3} =&-0.0904316832390810561,&
b_{1,3} =&-0.0676404519071381907,\\
b_{2,3} =&0.06759613017704596460,&
b_{3,3} =&0.02955525704293155274,\\
b_{6,3} =&-0.0000139180257516060,&
b_{0,4} =& 0,\\
b_{1,4} =& 0,&
b_{2,4} =&-0.0923364619367118592,\\
b_{3,4} =&-0.0169364939002081717,&
b_{6,4} =&-0.0000140086798182036.
\end{array}
\]

\

\paragraph{$k=6$ products}

With six products we can reconstruct the Taylor polynomial up to degree $n=22$ by applying the same strategy. We have also explored
different alternatives, considering decompositions based on the previous computation of low powers of the matrix, such as 
$A^2$, $A^3$, $A^4$, $A^8$, etc. to achieve degree $n=24$, but all our attempts have been in vain. Nevertheless, we should remark
that even if one could construct $T_{24}(A)$ with only six products, this would not lead to a significant advantage with respect to 
considering one scaling and squaring ($s=1$ in eq. (\ref{eses})) and using the previous decomposition for $T_{18}$.


In Table~\ref{tab.productsPS_NEW} we show the number of products required to evaluate $T_n$ following the PS rule and the new 
decomposition strategy. The improvement for $k \ge 3$ products is apparent. 
\begin{table}
\caption{\label{tab.productsPS_NEW}Minimal number of products to evaluate a polynomial of a given degree.}
\centering
		\begin{tabular}{lrrrrccc}
		\toprule
			Paterson-Stockmeyer\\
			Products & 1&2&3& 4&5&6  \\ 
			Degree   & 2&4&6& 9&12&16\\
			\midrule
			New decomposition\\
			Products & 1&2& 3&4 &5 &6\\ 
			Degree   & 2&4& 8&12&18&$22$\\
		\bottomrule
		\end{tabular}\\[5mm]
\end{table}

\section{Numerical performance and stability}
We next try to assess the performance of the various approximations  by comparing
the efficiency of the new procedure to compute the Taylor polynomial with the PS rule and also with diagonal Pad\'e
approximants, all in combination with  the scaling and squaring method for the matrix exponential.

All methods are based on an estimate of the matrix norm $\|A\|$, which is then used to choose the optimal order of the approximant together with the scaling parameter, based on the use of backward error bounds. 

\begin{figure}[t]
\centering
\caption{\label{fig.staircase_all}
Optimal orders and their respective cost vs. the norm of the matrix exponent. The numbers indicate the order of the approximants.
The diagonal lines show the asymptotic cost based on scaling and squaring and the corresponding orders are highlighted in bold-face.
}
\pgfplotsset{every axis plot/.append style={line width=1.0pt, mark size=3pt},
		tick label style={font=\footnotesize},
		every axis/.append style={%
		minor grid style={line width=.2pt,draw=gray!50},
		scale only axis, 
		}
	}
	\setlength\figurewidth{.8\textwidth}%
	\setlength\figureheight{.4\textwidth}
	\tikzsetnextfilename{figs/padepsnew}
\includegraphics{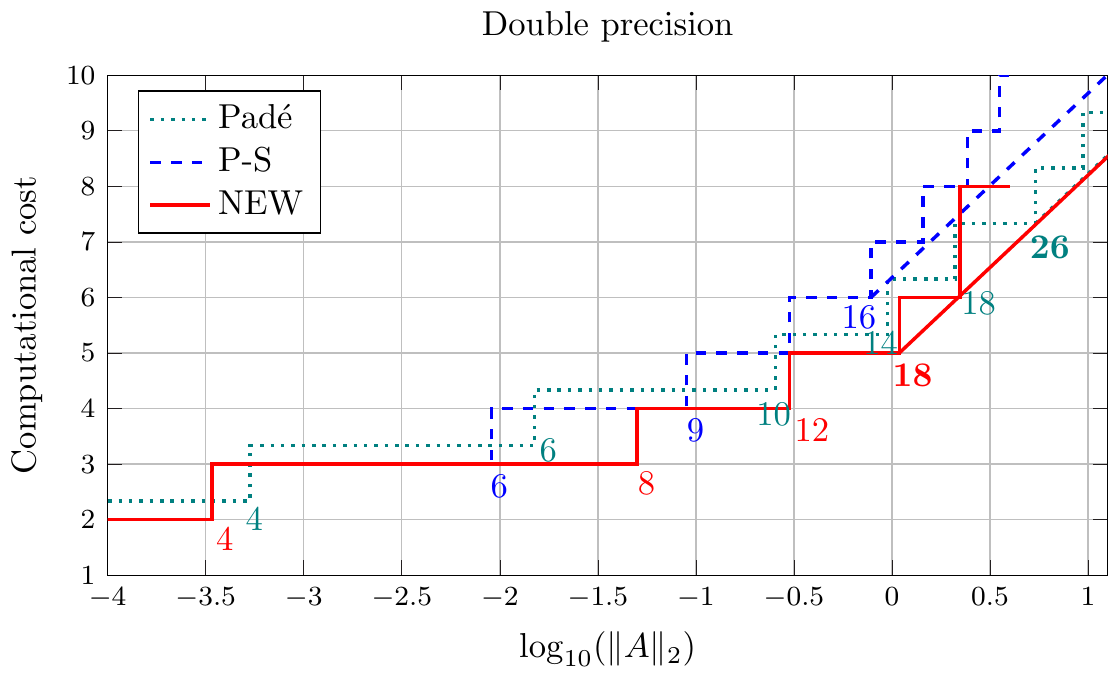}
\tikzsetnextfilename{figs/padepsnew_single}
\includegraphics{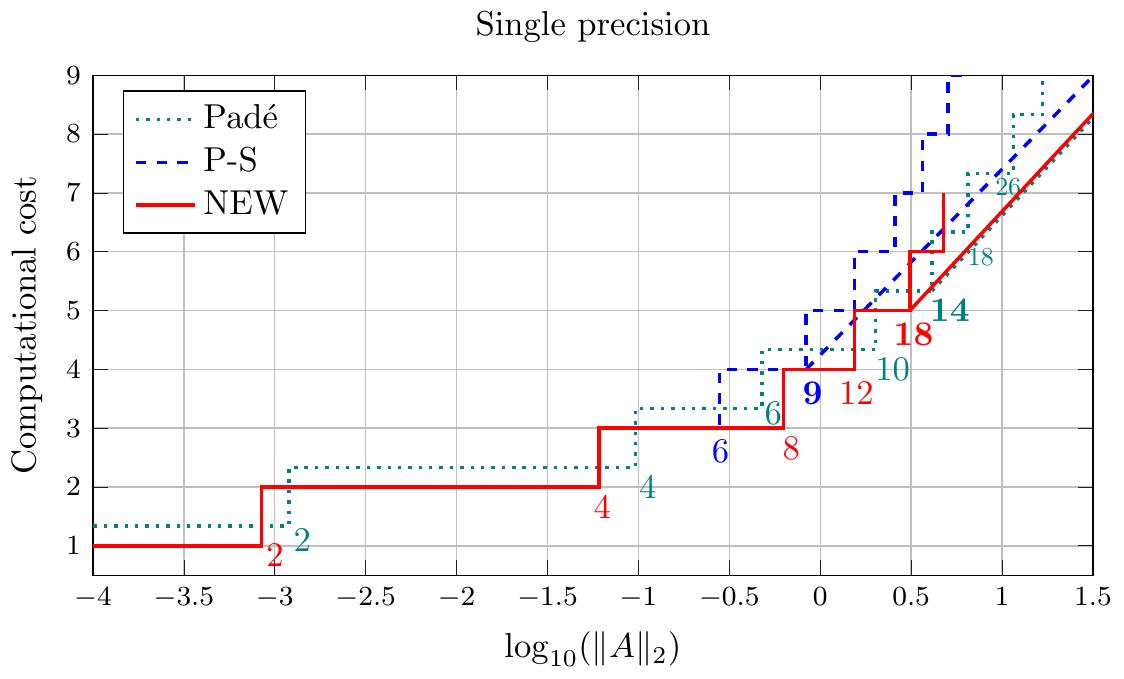}

\end{figure}

In Fig.~\ref{fig.staircase_all} we plot $\|A\|$ versus the cost measured as the number of matrix products necessary to evaluate 
Pad\'e and Taylor with the PS rule and the new implementation to approximate $\e^A$ in double precision (top) and single precision
(bottom). It is clear that for relatively small values of  $\|A\|$ the new approach to compute the Taylor approximant shows
 the best performance, whereas for higher values it has similar performance as Pad\'e. The diagonal lines in the graphs show the
 asymptotic cost based on scaling and squaring.
 
 %

\begin{figure}
\centering
\caption{\label{fig.stability}
Relative errors for a range of $30\times 30$ matrices computed with various algorithms as discussed in the text. Notice that all methods show similar (double) precision which is close to round-off accuracy. The error growth towards larger norms stems from the squaring error and affects all methods alike. The black horizontal line is the backward precision goal of the scaled exponential, cf. \eqref{eq.theta}. }
\pgfplotsset{every axis plot/.append style={line width=1.0pt, mark size=3pt},
		tick label style={font=\footnotesize},
		every axis/.append style={%
		minor grid style={line width=.2pt,draw=gray!50},
		scale only axis, 
		}
	}
	\setlength\figurewidth{.8\textwidth}%
	\setlength\figureheight{.4\textwidth}
	\tikzsetnextfilename{figs/error_double}
\includegraphics{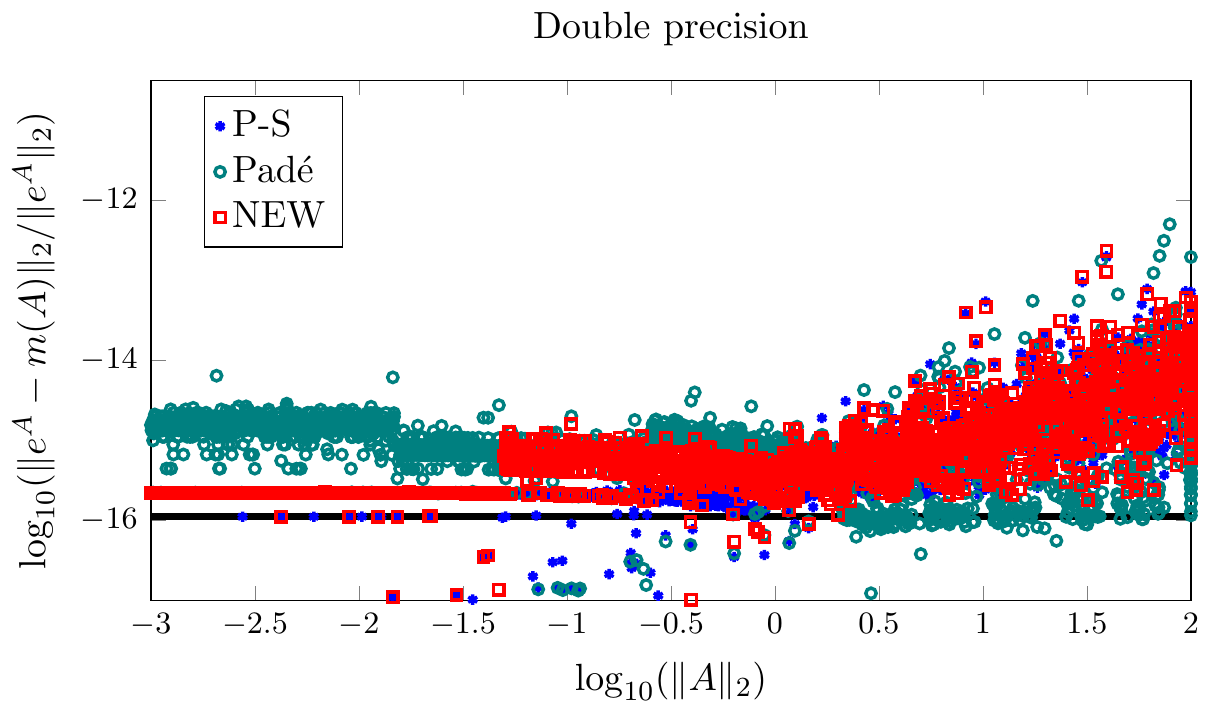}
	\tikzsetnextfilename{figs/error_single}
\includegraphics{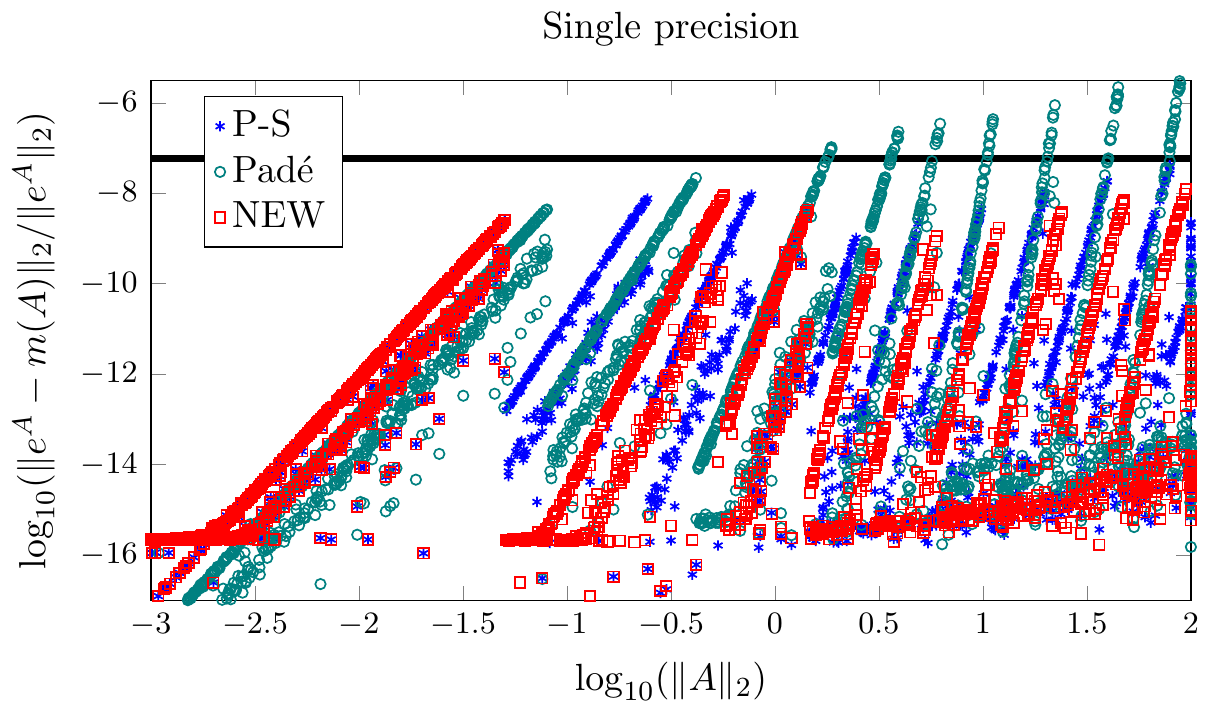}
\end{figure}

 To check the numerical stability of the proposed technique we have considered 46 different special matrices from the {\sc Matlab} matrix 
collection, as has been done in \cite{higham05tsa}, sampled at 1000 different norms, as well as 2000 random matrices which 
are scaled by random scalars to spread them over a range of norms. In  Fig.\ref{fig.stability} we plot the relative errors in the
computation of the corresponding matrix exponential vs. the norm of the matrix in double (top) and single precision (bottom) for 
the different approximants analysed here.
All matrices are of dimension $30\times30$, but the same experiments with matrices of larger dimension (up to $60 \times 60$) show
virtually identical results.

As a reference, we use a diagonal Pad\'e method of order 26 and the matrices are scaled and squared according to criterion \eqref{eq.theta} until the resulting norm is below the accuracy threshold from Table~\ref{tab.theta}.

We observe that the accuracy of our method lies well within the range that can be expected for Horner scheme evaluations and the Pad\'e algorithm from Table~\ref{tab.theta}.


\section{Concluding remarks}

It is important to remark that very often the evaluation of the exponential of a matrix constitutes an intermediate stage in the process
of solving differential equations by exponential integrators. This is the case, in particular, for equations of the form
\[
\begin{aligned}
	& X'=F(t,X)X, \\
	& X'=[F(t,X),X],  \\
	& X'=[[F(t),X],X], 
\end{aligned}
\]
where $X(t_0)=X_0\in\mathbb{R}^{n\times n}$, $t\in[t_0,T]$ and $[ \cdot , \cdot]$ denotes the usual matrix commutator.
 Among the advantages of exponential integrators applied to these equations are that they preserve
qualitative properties and their favourable long-time error propagation \cite{hockbruck10ein}. 

Take for example the linear equation $X' = F(t) X$. If $F(t)$ belongs to some particular matrix Lie algebra (e.g., $F$ is skew-symmetric),
then the solution $X(t)$ evolves in the corresponding Lie group (e.g., $X(t)$ is orthogonal for all $t$), and one would like that any
numerical approximation to $X$ obtained by applying an integrator share this important feature with the exact solution. If the integrator
is formulated in terms of matrix exponentials, such as with Magnus or Fer integrators \cite{blanes09tme}, this is achieved by construction. 

Given a time step, say $h$, one usually considers a numerical integrator of this class of order $h^p$, for values of $p$ usually 
between 2 and 6. Since the matrices involved in the schemes are typically proportional to $h$
and the truncated polynomial approximation, $T_n$, is an approximation to order ${\cal O}(h^{n+1})$, then one rarely needs to consider $n>18$.


In some cases one can take advantage of the very particular structure of the matrix $A$ and then design especially
tailored (and very often more efficient) methods for such problems \cite{celledoni00ate,celledoni01mft,iserles00lgm}.
Also, if one can find an additive decomposition $A=B+C$ such that $\|C\|$ is small and $B$ is easy to exponentiate, e.g., $\e^D$ is sparse and exactly solvable (or can be accurately and cheaply approximated numerically), and $B$ is a dense matrix, then more efficient methods can be found \cite{bader15tss,najfeld95dot}.
In addition, if $A$ is an upper or lower triangular matrix, it is shown in \cite{almohy09ans} that it is advantageous to exploit the fact that the diagonal elements of the exponential are exactly known. 
It is then more efficient to replace the diagonal elements obtained numerically by the exact solution before squaring the matrix. This technique can also be extended to the  first super (or sub-)diagonal elements.


In summary, we have derived a new algorithm to compute the matrix exponential by reducing the cost of the evaluation of 
its Taylor polynomial.
The method is as stable as the Pad\'e methods implemented in the {\sc Matlab}  function \verb#expm#  and 
is cheaper to compute than Pad\'e approximants for a wide range values for the matrix norm. For the convenience of the reader,
we provide a {\sc Matlab} implementation of the proposed scheme on our website \cite{web.algo}.

\subsection*{Acknowledgments}

This work has been funded by Ministerio de Econom\'{\i}a, Industria y Competitividad (Spain) through project MTM2016-77660-P (AEI/FE\-DER, UE).

\end{document}